\documentclass[]{article}
\usepackage{amsmath,amsfonts}
\usepackage{graphicx}
\usepackage[margin=2.5cm]{geometry}
\newcommand{\mrd}{\mathrm{d}}
\usepackage{hyperref}

\title{Histogram tomography}
\author{W.R.B Lionheart}

\begin{document}
 
\maketitle

\begin{abstract}
In many tomographic imaging problems the data consist of integrals along lines or curves. Increasingly we encounter ``rich tomography" problems where the quantity imaged is higher dimensional than a scalar per voxel, including vectors tensors and functions. The data can also be higher dimensional and in many cases consists of a one or two dimensional spectrum for each ray.  In many such cases the data contain not just integrals along rays but the distribution of values along the ray. If this is discretized into bins we can think of this as a histogram. In this paper we introduce the concept of ``histogram tomography''. For scalar problems with histogram data this holds the possibility of reconstruction with fewer rays. In vector and tensor problems it holds the promise of reconstruction of images that are in the null space of related integral transforms. For scalar histogram tomography problems we show how bins in the histogram correspond to reconstructing level sets of function, while moments of the distribution are the  x-ray transform of powers of the unknown function.  In the vector case we suggest a reconstruction procedure for potential components of the field. We demonstrate how the histogram longitudinal ray transform data can be extracted from Bragg edge neutron spectral data and hence, using moments, a non-linear system of partial differential equations derived for the strain tensor. In x-ray diffraction tomography of strain  the transverse ray transform can be deduced from the diffraction pattern the full histogram transverse ray transform cannot. We give an explicit example of distributions of strain along a line that produce the same diffraction pattern, and characterize the null space of the relevant transform. 
\end{abstract}

\section{Introduction}

Many tomographic imaging problems fall in to the broad category of {\em rich tomography} problems meaning that for each line our data are more than a single scalar. For example in many absorption tomography problems we might record a spectrum, a real valued function on an interval corresponding to the amount of radiation at each frequency or wavelength transmitted along that line. In other examples we might apply a narrow beam of radiation, x-rays, neutrons or electrons, and measure a diffraction pattern -- a real valued function of two variables. From this data we expect to be able to recover more than a single scalar image. We might expect a vector quantity such as velocity or magnetic field, a tensor such as strain, or the concentration of more than one chemical. Such techniques have been described as {\em Rich tomography}.  Alternatively we might expect to recover a single scalar image with fewer projections than would be necessary for traditional scalar measurements. In many cases acquiring a projection from raster scanning a beam is time consuming and reconstruction with fewer projection desirable, in other cases where rotation about an axis is limited it is essential. 

In particular numerous problems where the data are spectral allow us to infer the {\em distribution} of a scalar variable along each ray from the measurement of the spectrum for that ray, rather than simply its integral. In a practical setting this will often mean a histogram, in which the relative frequency of each value along the line that falls in to each bin of the histogram is recorded. In the limiting case as the size of the bins vanish this is the distribution, explained in the next section.

As a general rule standard tomographic problems involve a transport equation along the ray. In this case the values recorded are cumulative. As a ray travels from one voxel to the next the transport equation transforms the value input to one voxel to that output to the next. In the simplest case of attenuation tomography this results in Beer-Lambert law and the logarithm of the measurement at at the output is a line integral \cite{Natterer}. In other {\em non-Abelian} tomography problems the transport equation cannot be solved using an exponential and we have a non-linear integral operator along rays. Polarimetric neutron spin tomography provides an example \cite{sales2018three,desai2019polarimetric}. The problems we treat in this paper of a quite  different nature. As an archetypal example consider infrared absorption tomography\cite{An,Ma} where a certain chemical species absorbs a narrow band of infrared light, but the centre frequency of that band is a known monotonic function of temperature. Assuming the species is uniformly  distributed over the medium a portion of the spectrum around the shifted bands, when transformed to the temperature variable, directly gives us the distribution of temperatures. We know which temperatures occur along the ray, and how often they occur on the ray, but from one projection we do not know the location where these temperatures occur. Interestingly we do know the maximum and minimum temperatures immediately.

In the next section we review the definition of the distribution of a function with respect to a measure, and the concept of moments of the distribution. This has already been explained very clearly in the case of the Doppler tomography of vector fields by Andersson \cite{andersson2005doppler}, who used the idea of moments of the distribution. We go on to consider scalar histogram tomography problems  and what can be deduced from limited data. We then review concepts from Sharafutdinov \cite{sharafutdinov} of symmetric tensor fields and their ray transforms. We return to the Doppler moment transform and show that from the second moment we can reconstruct a  potential field directly. We then review the method of using the neutron transmission spectra near a Bragg edge and show that in principle the histogram of the strain in the ray direction can be recovered from this data. It was mentioned in \cite{LW} that as the linear elastic strain is potential the longitudinal ray transform derived from the Bragg edge yields no data on the interior strain. We show that the moment data of the distribution gives rise to non-linear partial differential equations in the strain and give explicitly the second moment in the two-dimensional case. We go on to discuss  x-ray strain tomography and the potential for extracting partial histogram data in this case, giving an explicit example of two strain distributions along a line resulting in the same diffraction pattern.

In this paper we have not specified the smoothness assumptions on functions and kept functional analysis to a minimum. This is an effort to make the paper accessible to the diverse applications communities that work with histogram tomography data. The details are easily filled in by anyone familiar with the necessary analysis for the Radon transform \cite{Natterer}, or the ray transforms of vector and tensor fields~\cite{sharafutdinov,boman2018stability}. 

\section{Distribution and moments}

\begin{figure}\begin{center}
		\includegraphics[scale=0.5]{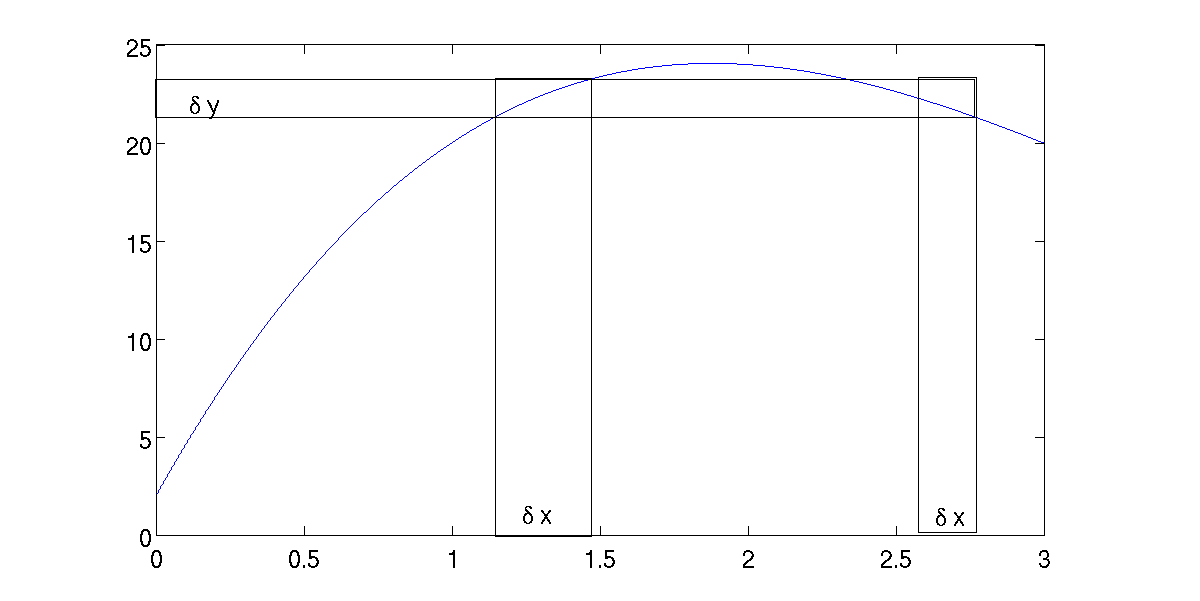}
	\end{center}
	\caption{\label{fig:distro-explain} The preimage under a function $f$ of an interval of length $\delta y$ in this case is two intervals with different lengths. In the limit as $\delta y \rightarrow 0$ each $\delta x = \delta y/|f'(x)|$. The distribution $\phi_f$ is the sum of such terms as given in (\ref{eq:phif}).}
\end{figure}

The idea of the probability distribution of a random variable is  familiar, but in this paper we will need the concept of the distribution of a function with respect to the standard measure on Euclidean space rather than a probability measure. The (Lebesgue) measure $\lambda$ on Euclidean space $\mathbb{R}^n$ is a function that assigns a measure of size to a set \cite{bogachev2007measure}, for example  total area ($n=2$), volume ($n=3$). For $n=1$ the measure of an interval $\lambda[a,b] = b-a$ is simply the length. 

\begin{figure}\begin{center}
		\includegraphics[scale=0.9]{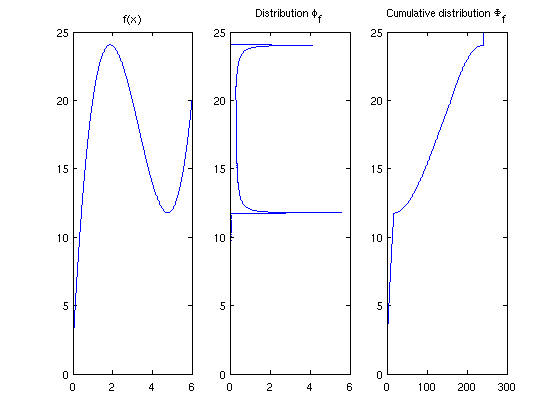}
	\end{center}
	\caption{\label{fig:distribution}Illustration of distribution and cumulative distribution of a smooth function $f(x)= (x -1)(x-3)(x-6)+20$. The distributions are shown with the frequency on the horizontal axes for easy comparison with the graph of $f$.  The plot is an approximation generated using the Matlab functions {\tt histc} and {\tt cumsum} using a discrete set of $x$ values.   Note that critical values of $f$ result in singularities of the distribution. In general the distribution gives how `frequently'  each value in the range occurs, but not where they occur in the  domain of $f$.}
\end{figure}
 Given  
a real valued  (measurable) function $f$ on a bounded measurable subset $\Omega \in \mathbb{R}^n$ we define the {\em push forward measure} $f_*\lambda$ on   $\mathbb{R}$ where for any interval $[a,b]\in \mathbb{R}$ $f_*\lambda [a,b] = \lambda(f^{-1}[a,b])$. Under mild hypotheses there is a (generalized) function $\phi_f$, the {\em distribution function} of $f$ such that
\begin{equation}
(f_*\lambda)[a,b] = \int\limits_a^b \phi_f(y) \,\mathrm{d}y.
\end{equation}

This definition is the same as the {\em probability distribution} of a random variable $f$ except that we do not assume a probability measure. 
We also define the cumulative distribution 
$ \Phi_f(y) = f_*\lambda (-\infty,y]$, 
the measure of the sub-levelset of $f$ up to $y$. Even though $\Phi$ need not be continuous we have $\Phi'_f = \phi_f$ in the sense of generalized functions. 
See figure (\ref{fig:distribution}) for an example of the distribution and cumulative distribution of a function. 

In the case where $f(x)=c$ a constant then clearly $\phi_f(y)= \delta(y-c) \lambda(\Omega)$. By contrast  for a differentiable function $f$ and (for simplicity dimension $n=1$) let $y\in \mathbb{R}$ be such that $f^{-1}(y)$ is finite, and $f'(x) \ne 0$ for all $x \in f^{-1}(y)$ (not a critical value), then
\begin{equation}\label{eq:phif}
\phi_f(y) = \sum\limits_{x \in f^{-1}(y)} \frac{1}{|f'(x)|}.
\end{equation}
See Fig \ref{fig:distro-explain} for an intuitive explanation.  Note that at points $y$, for which there is an $x$ with $y=f(x)$ and $|f'(x)|$ is small, $\phi_f(y)$ will be large.

As illustrated in Figure \ref{fig:distribution} isolated critical points of a smooth function result in singularities of the distribution. This means that one can in principle count distinct critical points and know their critical values from the distribution, and so approximately from a finely binned histogram. Of course one cannot infer the location of the critical points from the distribution.

The first moment of the distribution of a function $f$,
\begin{equation} \bar{f}=
\int_{\mathbb{R}} y \phi_f(y) \,\mathrm{d}y = \int_{\Omega} f(x) \,\mathrm{d}x ,
\end{equation}
is the integral of the function. If $f$ were a random variable and the measure normalized to be a probability measure, this would be the expected value. We also define the $k$-th moment of the distribution of $f$ for $k \in \mathbb{N}$ as
\begin{equation} \label{eq:mkf}
m_kf= \int_{\mathbb{R}} y^k \phi_f(y) \,\mathrm{d}y = \int_{\Omega} f(x)^k \,\mathrm{d}x 
. \end{equation}
For positive monotonic smooth functions (\ref{eq:mkf}) is simply change of variables for integrations as (\ref{eq:phif}) has only one term in the sum, while  \cite[Cor.  A.1.]{andersson2005doppler} proves it for general functions that can take negative values.

We will need to consider in some cases generalized functions, such as Dirac $\delta$-functions, which are also called distributions. To avoid confusion with the distribution of a function, in this paper we use the term generalized functions.

%


Let $L_{\theta,p}=\{x:x_1\cos \theta + x_2 \sin \theta =p\} =\{ x:x\cdot \Theta =p \}$, where $\Theta=(\cos \theta,\sin \theta)$, be a line in the plane.
We define the Radon transform on functions on the plane by
\begin{eqnarray}
\label{radonTransform}
Rf(\theta,p) = \int\limits_{-\infty}^\infty f(p \Theta +s \Theta^\perp) \,\mrd s
= \iint_{\mathbb{R}^2} f(x) \delta(\Theta\cdot x-p) \mrd x_1 \mrd x_1 
\end{eqnarray}  
where $\Theta^\perp= ( -\sin \theta, \cos \theta)$.  We generalize this to the {\em histogram Radon transform}
\begin{equation}
Hf (\theta,p,y) = \phi_{f|_{L_{\theta,p}}}(y)
\end{equation}
 the distribution of $f$ restricted to each line. Clearly we can recover the Radon transform by taking the first moment of the distribution with respect to $y$
\begin{equation}
Rf(\theta,p)  = \overline{Hf (\theta,p,\cdot)}.\end{equation}
We note also that 
\begin{equation}\label{momentradon}
m_k Hf(\theta,p,\cdot)= R (f^k) (\theta,p).
\end{equation}

  This shows that for the scalar case each of the moments of $Hf$ produces no more information than the Radon transform of the same power of $f$ and for a non-negative function $f$ exactly the same data, from invertability of the Radon transform. Indeed a non-negative bounded function is determined completely by the set of all of its  moments \cite{akhiezer1965classical} so the data are identical. It is interesting to note however that while fitting a function $f$ to its histogram tomography data $Hf$ is a non-linear problem, each of the problems (\ref{momentradon}) is linear.

One of the main  questions, for the scalar case in the plane, is to identify interesting subsets of lines for which $Hf$ determines $f$, without the reduction in stability associated with the limited data problem for the Radon transform.

We will need the Radon plane transform of functions of 3-space later. This is simply the integral over a plane normal to a unit vector $\Theta \in \mathbb{R}^3$ a distance $p$ from the origin so the notation is the same as for the integral over lines in the plane
\begin{equation}
Rf(\Theta,p) = \int_{x: \Theta\cdot x=p} f = \int_{\mathbb{R}^2} \delta(\Theta\cdot x-p) f(x) \,\mrd x
\end{equation}

\section{Recovery of level sets and discrete tomography}

A key insight that provides the connection between between the histogram tomography transform and the body of classical work on tomography is that the cumulative histogram gives the Radon transform of sub-level sets. To unpack that statement, the sub-level set of $f$ at $y$ is  $S_f(y)=\{ x: f(x)\le y\}$,and
$$(f|_{L_{\theta,p}})_* \lambda (-\infty,y] = \int\limits_{-\infty }^{y} Hf(\theta,p,y') \,\mathrm{d} y' = R\chi_{S_f(y)}(\theta,p)$$
is just the total length (measure) of the intersection of $L$ with the set which is also the Radon transform of the characteristic function $\chi_{S_f(y)}$. 

This leads us to the  {\em chord length problem}, which  is  to recover a bounded subset $S$ in the plane from this (exact) data for some set of lines. For convex sets $S$ Hammer's x-ray problem is to recover $S$ from $R\chi_{S}(\theta,s)$ for a specified set of $\theta$ and all $s$. For convex $S$ the solution to this problem is unique for four values of $\theta$ that are such that the cross-ratio of their slopes is a transcendental number \cite{GardnerMcMullen1980}. More practical examples of sets of four directions that are sufficient to determine a convex set are given in \cite{gardner1997discrete}. A reconstruction algorithm for convex polygon with $k$ sides from $k$ projections, possible with noisy data, is given in \cite{gardner2007solution}. See also  \cite{gardner2006geometric}. 

For non-convex sets less is known. But there is another approach. Local tomography methods use limited data for the Radon transform to recover the singular support of a function. In the case of  $\chi_{S_f(y)}$ the singular  support is exactly the contour $f^{-1}(y)$. As an example \cite{faridani1992local} shows that it is only necessary to know Radon data for all the lines passing through a neighbourhood of the singular support to recover it uniquely. This means for each contour one only needs the subset of the histogram tomography data for lines passing close to that contour $y$ value. 

As characteristic functions take only the values 0 or 1,  Radon transform inversion for these functions  falls into the approach of what is called by Herman and Kuba {\em discrete tomography}  \cite{herman2012discrete}, that is tomography where the range of the unknown function is a finite set. In the $\{0,1\}$ case this is also called {\em geometric tomography} by Gardner \cite{gardner2006geometric} with the emphasis mainly on convex bodies. As far as numerical algorithms for discrete tomography the Discrete Algebraic Reconstruction Technique (DART) is a popular variation of the standard iterative algorithm used in tomography adapted to functions with a binary range \cite{batenburg2011dart}.

As mentioned in Sec 2 the distribution of a function along a line has singularities corresponding to critical values, so if there are isolated critical points with distinct critical values these can be seen on the distribution. Suppose now we have a function $f$  on the plane with isolated critical points. Let $x_c$ be a critical point with $f(x_c)=y_c$  and  $\nabla f (x_c)=0$ then for any line through $x_c$, that is $L_{\theta,p}$ such that $\Theta\cdot x_c =p$ the distribution $Hf(\theta,p,\cdot)$ will have a singularity at $y_c$. This could be used to gain information about critical points from a limited data from the histogram Radon transform.

As well as discretization of the range of values of $f$ it is common to discretize the domain in to pixels or voxels. If the pixels are reduced to points and the rays still taken to have no thickness this results in tomography problems on finite subsets of integer lattices, also known as discrete tomography\footnote{In a personal communication Richard Gardner points out that the term ``discrete tomography'' was originally coined by Larry Shepp}
but in this case with  more of a  number theory emphasis. 
	  
	  In a more general setting the  function is replaced by a function on the vertices of a graph and a set of paths through edges take the role of rays. 
 In one case the values of the function can take restricted to a  discrete set, in this context refereed to as colours.  The ``rays'' are simply subsets of vertices. The number of each colours occurring in each ray is specified we get a {\em discrete tomography graph colouring problem}. Mainly the computational complexity of such problems is studied in \cite{bentz2008graph}, rather than specific algorithms. It is however notable that in a fully discrete setting this is the histogram tomography problem we describe with labels rather than numerical values. Such problems arise in scheduling of work rosters  and have not drawn much attention from those working on spectral tomography methods.

\section{Vector and tensor ray transforms}

An interesting extension of scalar tomography problems is the tomographic imaging of vector and tensor fields. Anisotropic material properties appear in a variety of physical contexts including isotropic materials subjected to deformation or flow, materials with a crystalline (or liquid crystal) structure, and layered or fibrous media where the structure is on a scale below the resolution of measurement (so appear anisotropic). Recovering a velocity field in a fluid from Doppler measurements gives an important example of vector tomography  \cite{schuster200820,sparr1995doppler}. The propagation of electromagnetic waves in a quasi-isotropic medium, one in which the anisotropy of the permittivity is small with small derivatives, was studied by \cite{kravtsov1968}. A detailed asymptotic analysis \cite{sharafutdinov} shows that this can be approximated by a system of ODEs along geodesics, which can in turn be approximated by and integral along geodesics.  Such problems arise in propagation of radio waves in the ionosphere as well as polarized light through birefringent media.
One application is the tomographic imaging of strain in glass or transparent polymers  \cite{aben2005photoelastic,tomlinson2006design,lionheart2009reconstruction}.  In crystalline or polycrystalline materials the strain tensor affects the diffraction of x-rays, neutrons and electrons and observation of the diffracted or transmitted rays holds the possibility for imaging the strain. For x-rays see \cite{LW}, electrons \cite{johnstone2017nanoscale} and a review of neutron diffraction  \cite{woracek2018diffraction}.

For a vector or tensor field there are a number of different tomography transforms we could define, the simplest would just be to take the scalar x-ray transform of each component.  In more interesting cases some projection of the vector or tensor is taken at each point along the ray, which depends on the direction of the ray,  and the components of this are integrated along the ray.

We will adopt a notation close to that of that used Sharafutdinov in his comprehensive book \cite{sharafutdinov}. This notation tends to be geometric and not specific to a particular application. First we introduce the simplest idea: just to contract the tensor with the direction vector of a ray resulting in a scalar.
Given a vector field $f$ the longitudinal ray transform (LRT) is defined as
\begin{equation}
If(x,\xi)= \int\limits_{-\infty}^\infty \xi \cdot f(x+s \xi)  \, \mathrm{d} s,
\end{equation}
while for a rank-2 symmetric tensor field $f$ the LRT is

\begin{equation}
If(x,\xi)= \int\limits_{-\infty}^\infty \xi \cdot f(x+s \xi) \cdot \xi \, \mathrm{d} s
\end{equation}
for $x,\xi \in \mathbb{R}^n, |\xi| =1$ where $\cdot$ denotes contraction. Note that parameterizing rays by $(x,\xi)$, through $x$ in the direction $\xi$,  has some redundancy as one can add any multiple of $\xi$ to $x$ and obtain the same ray. In three dimensional space it is easily seen that the space of rays is a four dimensional manifold. For example one can derive a coordinate chart by considering only $x$ in some plane with unit normal $\nu$ and any unit vector $\xi$ in the hemisphere with $\xi\cdot\nu>0$.

 The definition of the LRT is extended to rank-$k$ symmetric tensor fields in the same way with the contraction, with respect to $\xi$,  $k$ times. As we have contraction and  our data consists of only one scalar per ray, we have not more data than the scalar x-ray transform.  It is not surprising then that the LRT has a null space. A simple argument using the fundamental theorem of calculus shows that {\em potential}
tensor fields, those that are symmetrized derivatives of of a lower rank tensor field are in the null space. Or tensor fields here are assumed to vanish at infinity.

 In the case of vector fields potential has just the usual definition: $f=\nabla u$ for a scalar $u$. Potential rank-2 tensor fields are those that can be expressed as $f = (\nabla u + \nabla u^T)/2$ for some vector field $u$. In general, following \cite{sharafutdinov}, we define the operator $\mrd$ from rank-$k$ to rank-$(k+1)$ symmetric tensor fields formed by differentiation and symmetrization. 

A rank-$k$ symmetric tensor field $f$ is said to be {\em potential} if $f=\mrd u$ for some rank-($k-1$) field $u$.  Within the Schwartz class the nullspace of the LRT is shown to be exactly potential fields \cite{sharafutdinov}. The same is known for  tensor fields  just in  $L^2$ spaces  on a convex domain $\Omega \subset \mathbb{R}^3$ \cite{boman2018stability}. In this case a potential field is defined as one of the form $\mathrm{d}v$ but in this case with $v|_{\partial\Omega}=0$

For vector fields we have the Poincar\'{e} lemma: on a simply connected domain potential vector fields are precisely those with vanishing curl.

For rank-$k$ symmetric tensors  the Saint-Venant tensor $W(f)$ plays the same r\^{o}le \cite{sharafutdinov}.
This is a  rank-$2k$ tensor  with the property that $Wf=0$ if, and on a simply connected domain  only if, $f$ is potential. For $k=1$, $(Wf)_{ij}= f_{i,j}- f_{j,i}$ is  the tensor of skew symmetric derivatives of $f$, with subscripts after a comma denote differentiation. This is the same as the exterior derivative of a 1-form or, after  rearrangement, of the three independent off diagonal components in dimension  $n=3$, the curl. Using the  permutation symbol and the convention that sums are taken over repeated indices the $i$-th component of the curl is $\epsilon_{ijk}f_{j,k}$.    For $k=2$, $Wf$ is the familiar Saint-Venant compatibility tensor that is well known in elasticity. By a similar rearrangement of components, for the case $n=3,k=2$ the non-zero components of the rank-4 Saint-Venant tensor can be rearranged as a rank-2 Kr\"{o}ner tensor $Kf$ (defined in Section \ref{sec:secondmoment}). For the general case (which we should call the Georgievskii-Kr\"{o}ner tensor) see \cite{georgievskii2016}. In general a rank-$k$ symmetric tensor can be expressed as the sum $f = f_{\mathrm{sol}} + f_{\mathrm{pot}}$. Here  $f_{\mathrm{sol}}$ is called the solenoidal part with $\delta f_{\mathrm{sol}}=0$, where $\delta$ is the generalized divergence. The potential part is $f_{\mathrm{pot}} = \mrd u$ for some $u$.   Note $Wf= Wf_{\mathrm{sol}}$  \cite{sharafutdinov}. 
%
%
With imposition of suitable far field or boundary conditions the decomposition  becomes uniquely determined. See also \cite{boman2018stability} for the less smooth case.

The importance of this decomposition is that from LRT data the solenoidal part can be recovered uniquely while the potential part is unknown.
  For $n\ge 2$ there is an explicit reconstruction formula for $Wf$ (or $Kf$) from $If$ of filtered back projection type \cite{sharafutdinov}. This determines the solenoidal part of $f$, so the solution is unique up to the  null space of potential fields.  The reconstruction formula as given in \cite{sharafutdinov} assumes {\em complete} data, that is for the full four dimensional space of rays in the case $n=3$ this is typically too expensive to collect in practical situations. However if one considers rays in a family of parallel planes one can recover the solenoidal part of the restriction of $f$ to each plane from that data.

In electromagnetic problems where the waves are transverse we find that we need to project tensor fields on the plane orthogonal to the ray direction. 
Let $P_\xi$ be the projection of a symmetric second rank tensor field on to the plane perpendicular to $\xi$, then the transverse ray transform (TRT) is defined as  
\begin{equation}
Jf(x,\xi)= \int\limits_{-\infty}^\infty P_\xi f(x+s \xi)  \, \mathrm{d} s.
\end{equation}
We consider the important case of dimension $n=3$. For   a direction $\eta$ normal to $\xi$
\begin{equation} \eta \cdot Jf(x,\xi) \cdot \eta = \int\limits_{-\infty}^\infty \eta \cdot f(x+s \xi) \cdot \eta \, \mathrm{d} s
\end{equation}
so in any plane normal to $\eta$ this is simply the scalar Radon transform of the component $\eta \cdot f\cdot \eta $. This means there is a simple reconstruction for data on planes normal to six suitably chosen \cite{LW} directions $\eta$. 

Both LRT and TRT  have histogram tomography version. For the HLRT the data is the distribution of $ f(x+t \xi) \cdot \xi $ for vectors and $\xi \cdot f(x+t \xi) \cdot \xi $ for rank-2 tensor. By contrast   $P_\xi f(x+t \xi)$, along the ray $x+t \xi$, is a tensor field so the HTRT would consist of the {\em joint distribution} of the values along the ray of the components of the tensor normal to the ray direction. 

In the case of the HTRT, one special case is that we have the distribution of $\eta \cdot f\cdot \eta $ along lines on a plane normal to $\eta$. As this in the integral version we have reduced to the scalar histogram tomography problem for the normal  component on this plane.  As the TRT, data consists of three independent components for each ray (equivalent to a symmetric $2\times 2 $ matrix). The distribution of these values can be considered as a  joint distribution. But for the argument just outlined only a {\em marginal distribution}, the distribution of the normal component to the plane, is needed. 

\section{Doppler velocimetry}

As Sparr \cite{sparr1995doppler} and Schuster \cite{schuster200820} explain Doppler velocity tomography data are already understood as the distribution of velocity components along a line, in the direction of a line. The technique is used in laboratory and field studies to determine a velocity field from the frequency shift of a narrow acoustic beam.  In our terms, after suitable signal processing, this is the HLRT of the velocity field. It is also called the Doppler spectral transform. The first moment is typically used and of course this gives only the solenoidal part of the velocity leaving the potential part to be determined by other means. Andersson \cite{andersson2005doppler} shows that the higher moments determine a system of non-linear partial differential equations for the potential part.
In this section we give a simple procedure to recover  potential velocity field  from the  HLRT.

Let $f= f_{\mathrm{sol}}+f_{\mathrm{pot}}$ be the decomposition of the velocity field in to potential and solenoidal parts. 
 As we can recover the solenoidal part of the vector field from the first moment we can assume that $f_{\mathrm{sol}}$ is known, hence we know the distribution $\phi_{L,\mathrm{sol}}$ of  $f_{\mathrm{sol}}\cdot\xi$ along a line $L$ in direction $\xi$. 
 Unfortunately this does not lead to a method to recover the distribution of the potential part as we do not know where along the lines the values arise. In probability theory there is a large literature on what can be determined about the distribution of the sum of two dependent random variables where the marginal distribution is known but not the joint distribution. In general if the marginal distribution is known for a multivariate random variable the joint distribution is defined only up to a {\em Copular} (the joint cumulative distribution after a change of variables so marginal distributions are uniform)
 by Sklar's theorem \cite{Sklar}, and the set of possible joint distributions with given marginals is called a Fr\'{e}chet class. The Fr\'{e}chet-Hoeffding Theorem gives bounds on the joint distribution \cite{nelsen2006copulas}.  While this suggests an interesting avenue of future work we will now confine ourselves to the special case of reconstructing a field with zero solenoidal part.

 In this case we  have   $f=\mrd u$ for a scalar $u$. Notice the second moment is nothing but the LRT of $\mrd u \odot \mrd u$. From \cite{sharafutdinov} we know we can recover  the Kr\"{o}ner  tensor \cite{georgievskii2016}
\begin{equation}
K_{mn} = \epsilon_{mik}  \epsilon_{nj\ell} \left(u_{,i}u_{,j}\right)_{k\ell} =  \epsilon_{mik}  \epsilon_{nj\ell} u_{,ik}u_{,j\ell},
\end{equation}
noting that the third derivative drop out as each have a pair of indices over which a skew symmetric sum is taken.
Recall that for a square matrix $A$, $A \mathrm{Adj}\,A = \det A I$ where $\mathrm{Adj}$ is the adjugate matrix, the transpose of the matrix of cofactors, and $I$ the identity matrix. For the $3\times 3$ case
\begin{equation}\label{adjadja}
(\mathrm{Adj}\,A)_{mn} = \frac{1}{2}\epsilon_{mik}  \epsilon_{nj\ell} a_{ik} a_{j\ell}
\end{equation}
and a simple calculation shows $\mathrm{Adj}\,\mathrm{Adj}\,A = (\det A)A$, and $(\det A)^2= \det \mathrm{Adj}A$. So when $\det A \ne 0$ (\ref{adjadja})  shows that we can deduce $A$ from $\mathrm{Adj}\,A$ up to sign as
$$A= \pm (\det \mathrm{Adj} A)^{-1/2}\mathrm{Adj}\mathrm{Adj}A=
\pm (\det \mathrm{Adj} A)^{1/2}(\mathrm{Adj}A)^{-1}.
$$
 When the determinant is zero in the case  $\mathrm{rank}\, A =2$ the  columns of $\mathrm{Adj}\,A$ span the null space of $A$. In the case $\mathrm{rank} \,A=1$ or $0$,  $\mathrm{Adj}\, A=0.$  

We see 
\begin{equation}
K = 2 \mathrm{Adj}\,\mrd^2 u
\end{equation}where  $(\mrd^2u)_{ij}=u_{,ij}$ second derivative matrix. 
When $\det K \ne 0$ we can recover $\mrd^2 u$ up to sign, but when the determinant is zero $K$ does not uniquely determine $\mrd^2u$. Explicitly  
we have that
\begin{equation}\label{duK}
\mrd^2 u = \pm \sqrt{\left( \frac{1}{2}\det K\right)} K^{-1}
\end{equation}
is known except for points at which $\det K=0$. We need some additional data to determine boundary conditions and the sign of $u$. Indeed it was pointed out in \cite{sparr1995doppler} that for radial rotation invariant vector field $f$ the HLRT of $f$ and $-f$ are identical. 

Note that $\mathrm{trace}\, \mrd^2u  = \nabla^2u$, so we consider the Poisson equation
\begin{equation}
\nabla^2 u = \pm\sqrt{\left(\frac{1}{2}\det K\right)} \mathrm{trace}\,K^{-1}
\end{equation}
and given boundary data for $u$ (for example outside an object it is zero) we can solve uniquely for $u$ given $K$, when $\det K \ne 0$, up to a sign ambiguity.
%
 
We now have a constructive procedure for recovering $u$ from the second moment of the histogram LRT of the potential vector field $\mrd u$ up to sign.  Apply the  reconstruction formula of \cite{sharafutdinov} to this data to recover  $K$
then solve Poisson's equation.

The extension to general velocity fields, with non-zero solenoidal part, is not as simple and we know little more than was already presented by Andersson\cite{andersson2005doppler}. The Solenoidal part can be assumed known from the first moment data. The Kr\"{o}ner tensor of the second moment is a second order system of partial differential equations for the potential with the solenoidal part as coefficients. We hope to investigate uniqueness of solution in subsequent work.

\section{Neutron Bragg edge strain tomography}

\begin{figure}
			\begin{center}
				\includegraphics[scale=0.5]{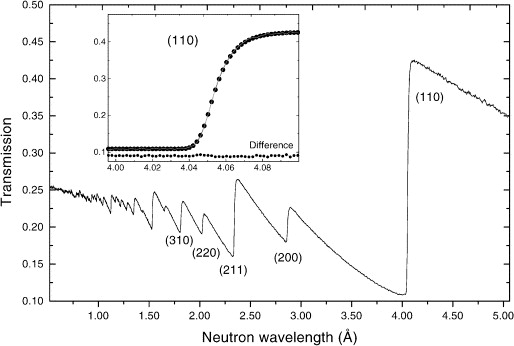}
			\end{center}
			\caption{\label{fig:Braggspectra}Neutron spectra from \cite{santisteban2002strain} showing the 110 edge on a magnified scale. The wavelength is given in Angstrom $= 10^{-10} m$. Note that the 211 edge is closer to the ideal behaviour of a step in an otherwise linear trend. }
		\end{figure}

Consider a sample of solid material that consists of identical crystals with a uniform random distribution of directions.  Each crystal has crystallographic planes denoted by Miller indices, which are triples of integers. A narrow beam of neutrons can be produced with a known distribution of kinetic energies. From a neutron scattering point of view the neutrons behave as waves with a wave length inversely proportional to  the energy, and this wavelength can be made to be close to the separation of the crystallographic planes. An incident wave is scattered  by parallel crystalographic planes with separation $d$ if the wavelength is a multiple of $2 d \sin \theta$. If one measures the neutrons that are transmitted, the scattered ones are lost. Looking at a transmission spectrum, such as Figure~\ref{fig:Braggspectra}, one sees a dramatic jump in the count of transmitted neutrons at a wavelength corresponding to $\theta=\pi/2$ for some crystallographic plane; this is called the Bragg edge. When the material is subjected to a linear elastic strain $\varepsilon$  the separation between planes with normal vector $\eta$ undergoes a relative change proportional to the strain $\eta\cdot\varepsilon \eta$ component.  See \cite{santisteban2002strain} for details and further references. If the strain was uniform along the path of the neutrons this would simply move the Bragg edge. 

\begin{figure}
	\begin{center}
		\includegraphics[scale=0.3]{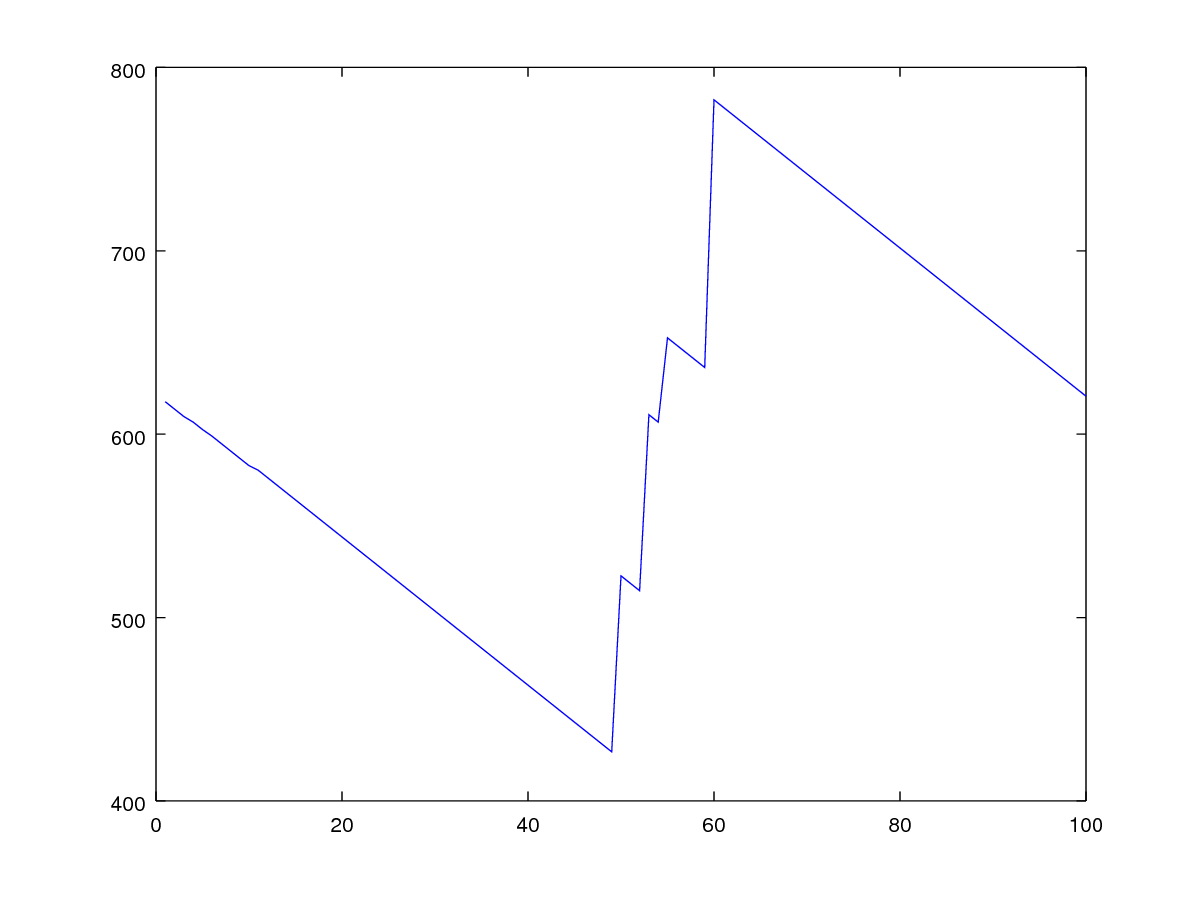}
	\end{center}

	\caption{\label{fig:braggshift} Shifted Bragg edges produce a saw tooth effect but the derivative is a histogram up to an added constant.}
\end{figure}

Under ideal conditions a Bragg edge can be modelled by linear function of wavelength multiplied by a Heaviside step function. Now consider a distribution of strains along the neutron path with a given histogram, the superposition of these Bragg edges is simply the sum of Heaviside functions, and hence the cumulative histogram, again multiplied by a linear function. On the enlarged plot of the Bragg edge corresponding to the Miller index 110 in Figure \ref{fig:Braggspectra} we see that this has the appearance of a sigmoid curve as expected of a cumulative distribution. To be precise we have to remove the linear component, which can be done by differentiation and then adjusting an added constant. See figure \ref{fig:braggshift} for a cartoon. We conclude that from this part of the spectrum we recover the HTLR of the strain tensor.  

Previously \cite{abbey2009} investigated the possibility of reconstructing strain from an average of the shift of a Bragg edge for each ray. The current author and Withers pointed out \cite{LW} that this amounts to the LRT of tensor field $\varepsilon = \mrd u$  in the interior of the object. This is not typically  potential in the sense of \cite{boman2018stability} where they consider a bounded convex  domain $\Omega$ and derive a decomposition of a symmetric rank $k$-tensor field in $L^2(\Omega)$  sum $f= f_{\mathrm{sol}} + \mrd v$ where $v|_{\partial \Omega}=0$ and $\delta f_{\mathrm{sol}}=0$, where derivatives are taken in sense of generalized functions.  They go on to prove that the solenoidal part $f_{\mathrm{sol}}$ can be recovered from  the LRT $If$, with a stability estimate using a Sobolev norm on the data. Instead we can use another result in \cite{boman2018stability} on fields in $L^2(\mathbb{R}^3)$ and write   $ \varepsilon = \chi_\Omega\mrd u $, extending by zero outside the domain. We still assume $\partial\Omega$ to be smooth but we know longer need it to be convex.  The decomposition means that $\chi_\Omega\mrd u =f_{\mathrm{sol}} + \mrd v$ with $v$ vanishing at infinity. The generalized function $f_{\mathrm{sol}}$ can be recovered uniquely and this determines some information about the boundary data for $\mrd u$.

This does at least provide some boundary data to solve the Lam\'{e} equations of linear elasticity, a second order elliptic system for $u$.
By standard uniqueness theory  $u$ is determined up an infinitesimal rigid motion by natural Neumann  data on the {\em stress} \cite{Knops1971}.  
Given LRT data plane by plane one could attempt to recover the strain from the measurements on rays in that plane using suitable additional a priori information.
One could represent the Lam\'{e} system by finite element or finite  difference methods provided the Lam\'{e} coefficients are known. Together with  a discetization of the LRT and a priori information an overdetermined linear system could be constructed and the least squares solution calculated.  This type of approach was taken for axisymmetric objects by Gregg {\it et al} \cite{gregg2017}, in the general case for synthetic data by Wensrich et al \cite{wensrich2016bragg}, and experimentally for a general case by Henriks {\it et al} \cite{hendriks20172}. 

%
%
%
%
%
%
%
		
The approach outlined above has some weaknesses. It is only known to work for convex objects, it uses only a summary (the first moment) of the distributional data, it requires the Lam\'{e} coefficients to be known,  and it requires the solution of a system of PDEs.

 Assuming now that the strain is purely potential, so that the the strain at the boundary of the object is zero, we now consider how we might use the the HLRT data contained in the detailed cumulative histogram in each neutron transmission spectrum. This supposes that we have enough neutron flux in the energy range around our chosen Bragg edge to have sufficient counts in each bin of energy to get some of the higher moments. So for some $k$ we have the $k$-th moment  and we see, analogously to the scalar case,  this is just  the  LRT of the symmetric tensor $k$-th power $\mrd u \odot \cdots \odot \mrd u$= $\mrd u^{\odot k}$. These are in general not potential so the moments will be non-zero.

So we can recover the Saint-Venant (or Kr\"{o}ner) tensor of these powers. As outlined in the next section this results in a nonlinear system of partial differential equations for each $k$. Unfortunately we know of no trick like that employed for vector fields that reduces these to a linear system. Essentially we are in the same situation described by \cite{andersson2005doppler}  for the vector case before this paper. While solving a nonlinear system is undesirable, and we expect more computationally demanding than solving the linear Lam\'{e} system, it does offer a new possibility at reasonable cost. The results of  \cite{hendriks20172} could be verified by computing these Saint-Venant tensors for several values of $k$ and comparing with the symmetric powers of the reconstructed strains. As there are few options for measuring the strain independently this may offer a valuable independent check on the results and hence the validity of the assumptions. The details of the calculations are in the next section.

\section{Second moment in the rank 2 case}\label{sec:secondmoment}

In general the Georgievskii-Kr\"{o}ner tensor for a rank $m$ tensor $f$ in dimension $n$ has rank $m(n-2)$. In the case $m=4,n=3$ it is given by
\begin{equation}
K_{i_{11}i_{21}i_{31}i_{41}} = \epsilon_{i_{11},i_{12},i_{13}} \epsilon_{i_{21},i_{22},i_{23}} \epsilon_{i_{31},i_{32},i_{33}} \epsilon_{i_{41},i_{42},i_{43}} f_{i_{13} i_{22} i_{33} i_{42},  i_{12} i_{23} i_{32}  i_{43}  }
\end{equation}
where we have followed Georgievskii \cite{georgievskii2016} in labelling indices with a pair of subscripts. Now consider the case $f = \mrd u \odot \mrd u$ where $u$ is a vector field 
\begin{equation}
f_{ijkl}  = u_{  ( i, j } u_{   k,l )}
\end{equation}
with round brackets indicating symmetrization over the indices.
In this case $K$ is a fully symmetric rank 4 tensor. In general the terms result from differentiating a product of derivatives of the $u_i$ four times. From the product rule we would expect terms including fifth derivatives of the components of $u$ but the skew symmetrizing over pairs of indices that appear as derivatives removes these terms.  In practice the components of $K$ consist of sums of pairs of third partial derivatives, and products of 4th and 2nd partial derivatives of $u$. For dimension $2$, or restriction to a plane, $K$ is a scalar consisting of a fourth order non-linear PDE for $u$. Unless we have additional information about $u$ this is insufficient to determine $u$ uniquely. However in the three dimensional case data for three or more families of planes normal to three independent vectors has the possibility of yielding a unique solution. Unlike the case of scalar $u$ we do not know a neat transformation that eliminates the non-linearity.


%
%

\section{X-ray diffraction strain tomography}

In the neutron Bragg edge case we ignored the scattered radiation, but it is possible instead to measure scattered particles/waves. In the Bragg edge case it was possible to use a  parallel beam of neutrons, as only the transmitted radiation was recorded.  But to capture a diffraction pattern for the material along a ray the simplest approach is to use a single narrow  beam, or possibly a small number of beams far enough apart that the diffraction patterns have large regions in which they do not overlap. In \cite{LW} we considered the case of a monochromatic x-ray beam and a polycrystalline material. We showed that by taking certain radial moments of the diffraction pattern we could recover the TRT of the strain. In this section we consider what histogram data is available form the diffraction pattern in this case. Intuitively the diffraction pattern is a function of two variables and yet the transverse strain has three independent components so we would not expect to recover the full HTRT of the strain.

Consider one specific ray in the direction $\xi$ and the projection of the strain in the plane normal to $\xi$ a  $2\times 2$ matrix valued function $A(s)$, as in \cite{LW}. The diffraction pattern for the voxel at $s$ is an ellipse determined by the matrix $A(s)$, and the total diffraction pattern the sum of these
\begin{equation}\label{integralforg}
g(y) =\int_{A:y^TAy=1} \phi(A)\,\mrd a_{11}\mrd a_{12}\mrd a_{22}/\sqrt{2}
\end{equation}
where $y$ is the coordinate in the plane of the screen where the diffraction pattern is formed and $\phi$ is the density function, which is a generalized function supported on the curve $A(s)$.  Consider a point $y=r(\cos \theta, \sin \theta)$ on the screen. The range of integration in (\ref{integralforg}) is then the points $(a_{11},a_{12},a_{22})$ satisfying the linear equation
\begin{equation}
\cos^2 \theta \, a_{11}  + 2 \sin \theta \cos \theta \, a_{12} +\sin^2 \theta \,a_{22} =r^{-2}
\end{equation}
The Frobenius norm $||A||_F^2 = a_{11}^2 + 2 a_{12}^2 + a_{22}^2$, so this agrees with the Euclidean norm on $\mathbb{R}^3$ if we use the coordinates $(a_{11}, \sqrt{2} a_{12}, a_{22})$. Let $n_\theta= (\cos^2 \theta,  \sqrt{2} \sin \theta \cos \theta ,\sin^2 \theta)$.  Notice $|n_\theta|=1$ and we see (\ref{integralforg}) is exactly the Radon plane transform 
\begin{equation}\label{RPT}
g(y) = R\phi(n_\theta,r^{-2})
\end{equation}
Note however that the family of planes is given by the two parameters $r,\theta$  where as the full set of planes in three dimensional space is three dimensional, so we do not expect to be able to recover a general $\phi$ from this. Of course the distribution $\phi$ is not completely general but is a generalized function supported only on a curve, so one might think the situation is somewhat better. Using the Fourier slice theorem for the ray transform we see that the data $g$ determines $\hat{\phi}$ on a cone. Let $\alpha=(\alpha_{11},\sqrt{2}\alpha_{12},\alpha_{22})$ be the Fourier transform variables then the components determined  by $g(y)$ are $\widehat{R\phi}(\alpha/|\alpha|,|\alpha|)=\hat{\phi}(\alpha)$ with $\alpha$ on the cone $\alpha_{11}\alpha_{22} = 4 \alpha_{12}^2$. We see that if the difference between two distributions has a Fourier transform vanishing on this cone then they will produce the same diffraction pattern $g$.

To make an explicit example we consider biaxial strain where the principle axes of the strain rotates along the ray though half a turn. The diffraction pattern $g$ will be supported and positive on an annulus with radii corresponding to the principle strains, and rotationally symmetric. Now consider a uniaxial strain in the transverse plane, in this case the bulk strain ranges between the principle strains of the previous case. This will also be positive and supported on the same annulus. To make the diffraction patterns the same we have to calculate the value of $g$ in the first case as a function of $r$  and arrange the bulk strain in the second example to produce the same answer.

Our case is 
\begin{equation}
A_1(s) = 
\left(\begin{array}{rr}
\cos s &-\sin s \\
\sin s & \cos s \\
\end{array}  \right)
\left(\begin{array}{rr}
2 & 0 \\
0 & 1 \\
\end{array}  \right)
\left(\begin{array}{rr}
\cos s &\sin s \\
-\sin s & \cos s \\
\end{array}  \right)  
=\left(\begin{array}{rr}
1 + \cos^2 s &\sin s \cos s\\
\sin s \cos s & 1+\sin^2 s \\
\end{array}  \right)
\end{equation}
Changing coordinates to $u=a_{11}+a_{22}, v=a_{11}-a_{22}, w = 2 a_{12}$ we see that $A_1$, in these coordinates, is $u=3, v= \cos 2s, w= \sin 2s$, that is a circle in the $u=3$ plane for $s$ ranging from $0$ to $\pi$. Calculation of $g$ is now the integral over the plane $n_\theta \cdot (a_{11},a_{12}/2,a_{22}) = r^{-2}$ of the generalized function $\delta(u-3)\delta( \sqrt{v^2+w^2}-1)$. In the new coordinates $U=(u,v,w)$ the plane under consideration in (\ref{RPT}) is 
\begin{equation}
\frac{1}{\sqrt{2}}( 1, \cos 2 \theta, \sin 2 \theta)\cdot U= 2 r^{-2}
\end{equation}
where the vector on the left has unit length.  This intersects the $u=3$ plane in a family of lines a distance $2\sqrt{2} r^{-2} -3$ from the $u$. The Radon transform of the unit density on the unit circle in $(v,w)$ space (see Appendix A)
\begin{equation}
R\delta(\sqrt{v^2 +w^2}-1)(p,\phi)= \chi_{[-1,1]}(p) \frac{2}{\sqrt{1-p^2}}
\end{equation}
To evaluate $g$ we need to take account of the Jacobian for the change of variables to $(u,v,w)$ and the direction cosine of $1/\sqrt{2}$ of the plane to the $u=3$ plane giving
\begin{equation}
g_1(r)= \frac{1}{\sqrt{2}} \chi_{[r_0,r_1] }(r)   \frac{2}{\sqrt{1 - \left(2\sqrt{2} r^{-2} -3 \right)^2 }}
\end{equation}
where $r_0 = 2/\sqrt{2+3\sqrt{2}}$, 
$r_1 = 2/\sqrt{4+3\sqrt{2}}$. We simply have to construct an $A_2(s)$ with a density supported on the $u$ axis that produces the same $g$. Each plane in (\ref{RPT}) intersects the $u$ axis at an angle $\pi/4$, at the point $U=(2\sqrt{2}r^{-2},0,0)$ so we just assign the density $\phi_2 = g_1(u^{1/2} 2^{-3/4})\delta(v)\delta(w)$. To find an example of an $A_2(s)$ with this density we simply form the cumulative density in the $u$ direction and use this as the bulk strain in the transverse direction, with a rescaling so that $s$ takes values in the same interval as for $A_1$.

One can  repeat this construction for any  pairs of curves with the property that any plane with normal in the cone passes through one point on each curve an by choice of the density function make the diffraction pattern the same. 

 In general the support of the joint distribution of a smooth transverse strain along a line will be an immersed curve in three dimensional space. When the planes defined above intersect this curve non-transversely a singularity will result in the marginal distribution that is given by any of the plane transforms in  the data. In the case detailed above a singularity occurs when the plane is tangent to the circle for the $A_1$ case and at the ends of the range on the $A_2$ case. 

Clearly we can recover the first moment and hence the TRT data as detailed in \cite{LW}. We can also recover the distribution of the normal component to a plane from the diffraction data for that plane. To see this suppose $a_{11}$ is the normal component. The directions $(1,0,0)$ satisfies the cone condition in the frequency domain, and so the Fourier transform along this direction determines the marginal distribution for $a_{11}$. Thus all results for the histogram Radon transform can be used plane by plane with this data. Similarly we know the marginal distribution of $a_{22}$. It is not yet clear how to use the additional data in an efficient way. In particular we would prefer a reconstruction method with fewer rotation axes than the three given in \cite{DesaiLionheart} for a general strain, or two for a potential (elastic) strain.

\begin{figure}
\centering
\includegraphics[width=0.4\textwidth,clip]{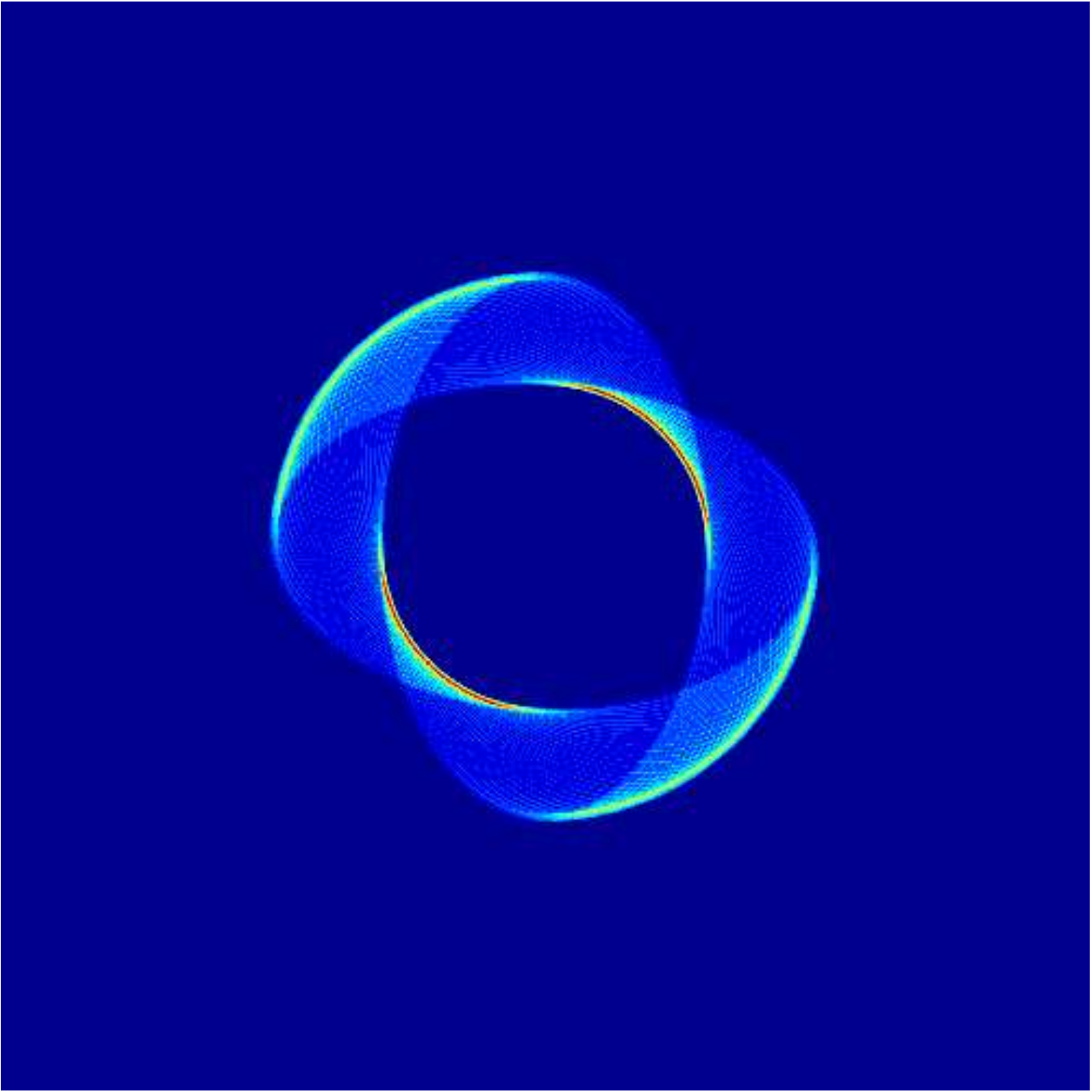}
\includegraphics[width=0.4\textwidth,clip]{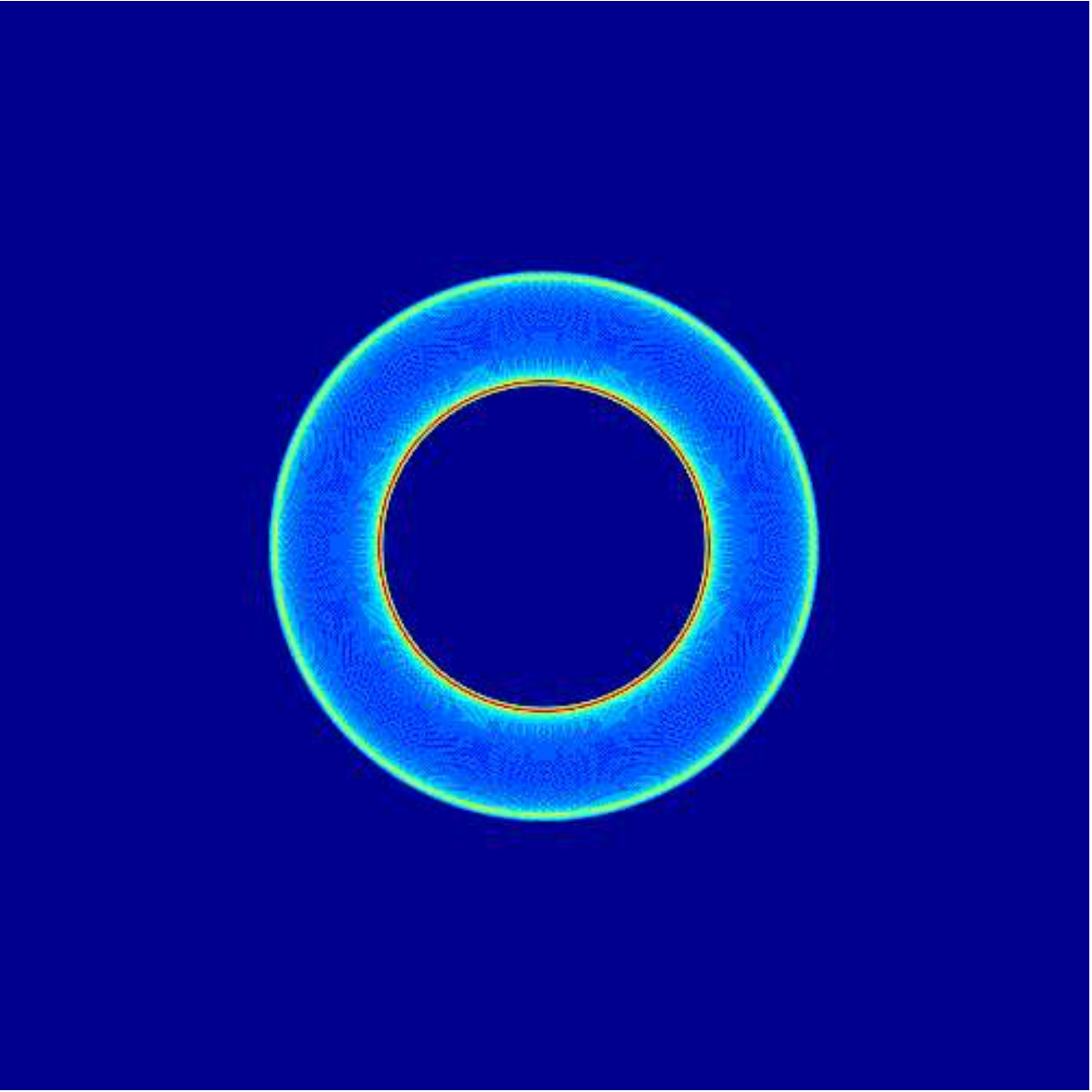}
\caption{Superposition of uniform density ellipses, left rotated through a quarter turn, right through a half turn.} 
\label{fig:ellipse_partial}
\end{figure}

\section{Conclusions and further work}

We have exhibited a wide range of tomographic  problems in which a distribution can be measured for each ray. While in general fitting a distribution (or in the discrete case a histogram) of the data is generally a non-linear problem we saw that in the case of a scalar unknown both the moment data and the bins in the cumulative histogram reduce our problem  to known linear tomographic problems.   Further work is needed on practical useful algorithms for scalar histogram tomography. For near infrared spectral tomography of chemical species, attention needs to be directed to the recovery of multiple parameters including the density of the chemical species and the pressure as well as temperature.

For a vector field, we were able to give for the first time an explicit reconstruction for the  potential case, using the second moment. The general case is an obvious candidate for further study, and especially the problem of limited data.

 In the case of Bragg edge strain tomography we showed how the histogram could be extracted from spectral data with a high resolution in energy near an edge. We were not so lucky with the theory as in the vector potential case, as we were only able to find a non-linear system of PDEs for the unknown displacement field derived from moment data. It still remains to be seen in practice what distributional data can be recovered from neutron spectra.
  
 In x-ray diffraction tomography of strain we were able to give an example of strain distributions resulting in the same diffraction pattern for one ray, so that in this case we do not recover the histogram transverse ray transform. We have yet to see experimental data demonstrating what can be recovered from x-ray diffraction data.
  
 We did not treat in this paper the case of strain tomography from  diffraction data from a single crystal. In contrast to the polycrystalline case the unstrained diffraction pattern for a pencil beam is a two dimensional lattice of small diffraction spots, and the lattice is distorted by a strain in the transverse direction. One might expect under ideal conditions and small strain to see each small spot in the lattice spread to a density function when the strain varies along the beam. The mean of these distributions gives the average displacement of the spot in reciprocal space. From this one would expect to be able to derive the TRT data of the strain. As the strain is potential, rotations about two axes are sufficient for reconstruction \cite{DesaiLionheart}. A practical case where this problem arises is electron strain measurement in silicon for electronic devices. The measurement technique is already established \cite{johnstone2017nanoscale}, however it is unlikely that measurements can be taken for a complete rotation in a sample of practical significance. This leads again to the consideration of reconstruction fro ma limited range of projections. Like the polycrystalline case, this technique would not yield  the full HTRT. Each spot gives a marginal distribution of a vector, but it is not clear what information about the joint distribution of the transverse strain tensor can be recovered.
  
  In many practical cases it is likely that algorithms that fit using all available data to the non-linear problem are better than simple explicit reconstruction formulae. Nevertheless theoretical methods at least give a useful insight in to the sufficiency or insufficiency of	 specific sets of measurements. 
  
  In practical problems there is always the response of the instrumentation and blurring effects from other physical phenomena to consider. One can generally expect the distribution that is measured to be convolved with a point spread function and some calibration and deconvolution is required to give a good estimate of a histogram. As histogram tomography is investigated experimentally this is one of the challenges that must me overcome. In many cases the first step will be to validate and calibrate the measured histogram from a sample with a known distribution.
  
  The more complex cases have forward problems involving restriction, distributions, linear projections and integrals. There are likely also to be problems where a weighted integral or histogram is taken, for example when scattered radiation is measured but it is also attenuated. Further work could look at a general theoretical framework that includes many practical cases. This might be expected to draw further on work from probability theory on recovering a joint distribution from marginal and conditional data. 
  
  We  look forward to practical realisations of these techniques as well as further theoretical developments and hope this paper will inspire such work.
  
\section*{Acknowledgements}
The author would like to thank the Royal Society for the Wolfson Research Merit Award that supports this work. Dr Nick Polydorides was very helpful in describing to the author in detail infrared spectral tomography, and thanks are due also to Prof Hugh McCann for organizing a workshop on Chemical Species Tomography in our MIRAN series funded by EPSRC (grant number EP/K00428X/1). Thanks to Prof Philip Withers for the introduction to x-ray and neuron strain measurement techniques; to Prof Alexander Korsunsky for discussion of Bragg edge neutron strain measurement; Prof Paul Midgley, Dr Carola-Bibiane Sch\"{o}nlieb  and colleagues at Cambridge for explaining the method of electron strain measurement; and Dr Sean Holman for helpful discussions on the mathematical theory. 
 
\appendix

\section{Radon transforms of generalized functions supported on curves}

The Radon transform is defined on generalized functions and the results are perhaps not obvious. Two methods for calculation are via the Fourier transform, and the Fourier slice theorem and using the back projection operator which is the formal adjoint of $R$ applied to a test function.

To calculate the plane transform of $f(x) =\delta(x_1)\delta(x_2)$  we have $\hat{f}(\xi)= \delta(\xi_3)$. Now $\widehat{Rf}(\Theta,\sigma) = \hat{f}(\sigma \Theta) = \delta (\sigma\Theta_3) = \frac{1}{|\Theta_3|}\delta(\sigma)$ and so $Rf(\Theta,s) = \frac{1}{|\Theta_3|}$

We now turn our attention to te generalized function supported on the unit circle $f(x) = \delta(|x|-1)$. In polar coordinates $\delta(r-1)$ so taking the Fourier transform $\hat{f}(\xi) = \int_0^{\infty}  \delta(r-1) J_0({|\xi|}r)r \mrd r  = J_0(|\xi|)$. Taking the inverse Fourier transform of the Bessel function $J_0$ we see
\begin{equation}
Rf(\theta,s) = \chi_{[-1,1]}(s) \frac{2}{\sqrt{1-s^2}}.
\end{equation}

\bibliographystyle{plain}
\bibliography{histo}

\end{document}